\documentclass[12pt]{article}
\usepackage{amsmath,amssymb}

\makeatletter
\def\eqnarray{\stepcounter{equation}\let\@currentlabel=\theequation
\global\@eqnswtrue
\tabskip\@centering\let\\=\@eqncr
$$\halign to \displaywidth\bgroup\hfil\global\@eqcnt\z@
  $\displaystyle\tabskip\z@{##}$&\global\@eqcnt\@ne
  \hfil$\displaystyle{{}##{}}$\hfil
  &\global\@eqcnt\tw@ $\displaystyle{##}$\hfil
  \tabskip\@centering&\llap{##}\tabskip\z@\cr}

\def\endeqnarray{\@@eqncr\egroup
      \global\advance\c@equation\m@ne$$\global\@ignoretrue}
\makeatother

\begin{document}
\bibliographystyle{tom}

\newtheorem{lemma}{Lemma}[section]
\newtheorem{thm}[lemma]{Theorem}
\newtheorem{cor}[lemma]{Corollary}
\newtheorem{voorb}[lemma]{Example}
\newtheorem{rem}[lemma]{Remark}
\newtheorem{prop}[lemma]{Proposition}
\newtheorem{stat}[lemma]{{\hspace{-5pt}}}

\newenvironment{remarkn}{\begin{rem} \rm}{\end{rem}}
\newenvironment{exam}{\begin{voorb} \rm}{\end{voorb}}

\newcommand{\ca}{{\cal A}}
\newcommand{\cd}{{\cal D}}
\newcommand{\ch}{{\cal H}}
\newcommand{\cs}{{\cal S}}
\newcommand{\cf}{{\cal F}}
\newcommand{\ra}{\rightarrow}
\newcommand{\spann}{\mathop{\rm span}}

\newcounter{teller}
\renewcommand{\theteller}{\Roman{teller}}
\newenvironment{tabel}{\begin{list}%
{\rm \bf \Roman{teller}.\hfill}{\usecounter{teller} \leftmargin=1.1cm
\labelwidth=1.1cm \labelsep=0cm \parsep=0cm}
                      }{\end{list}}

\newcommand{\Ni}{{\bf N}}
\newcommand{\Ri}{{\bf R}}
\newcommand{\Ci}{{\bf C}}
\newcommand{\Zi}{{\bf Z}}
\newcommand{\Ti}{{\bf T}}

\newcommand{\proof}{\mbox{\bf Proof} \hspace{5pt}}
\newcommand{\remark}{\mbox{\bf Remark} \hspace{5pt}}
\newcommand{\ruimte}{\vskip10.0pt plus 4.0pt minus 6.0pt}

\newcommand{\RRe}{\mathop{\rm Re}}
\newcommand{\ad}{\mathop{\rm ad}}
\newcommand{\Ad}{\mathop{\rm Ad}}
\newcommand{\Sp}{\mathop{\rm Sp}}
\newcommand{\supp}{\mathop{\rm supp}}
\newcommand{\Ker}{\mathop{\rm Ker}}
\newcommand{\AK}{}
\newcommand{\one}{1\hspace{-4.5pt}1}

\hyphenation{groups}
\hyphenation{unitary}

\thispagestyle{empty}

\begin{center}
{\Large\bf Approximately inner derivations} \\[5mm]
Ola Bratteli \\
{\small Department of Mathematics, University of Oslo}\\
{\small
Blindern, P.O.Box 1053, N-0316, Norway}
\smallskip\\ Akitaka Kishimoto \\ {\small
Department of Mathematics, Hokkaido University, Sapporo 060-0810,
Japan}\\ and\\
 Derek W. Robinson\\
{\small Centre for Mathematics and its Applications, Australian
National University}\\
{\small Canberra, ACT 0200, Australia}\\
{\small  March  2007}
\end{center}

\vspace{5mm}

\begin{center}
{\bf Abstract}
\end{center}

\begin{list}{}{\leftmargin=1.8cm \rightmargin=1.8cm \listparindent=10mm
   \parsep=0pt}
\item
Let $\alpha$ be an approximately inner flow on a $C^*$ algebra $A$ with generator
$\delta$ and let $\delta_n$ denote the bounded generators of the approximating flows $\alpha^{(n)}$.
We analyze the structure of the set 
\[
  \cd=\{x\in D(\delta): \lim_{n\rightarrow\infty}\delta_n(x)=\delta(x)\}
  \]
  of pointwise convergence of the generators.
  In particular we examine the relationship of $\cd$ and various  cores related to spectral subspaces.
\end{list}

\section{Introduction}

The theory of flows  on operator algebras has been largely motivated by
models of quantum statistical mechanics in which the flow  $\alpha$ is usually constructed as the
limit of a sequence $(\alpha^{(n)})$ of local flows.
The latter are typically given by inner one-parameter automorphism groups of the algebra $A$ 
and the flow is correspondingly called approximately inner.
In particular this is the situation for models of quantum spin systems in which the algebra of observables
$A$ is a UHF-algebra (see, for example, \cite{BR2}, Chapter 6).
An early result of  Sakai established \cite{S76}  that if $\delta$ is the generator of a flow $\alpha$ on a UHF-algebra then there exists an increasing sequence $(A_n)$ of finite dimensional
$C^*$-subalgebras $A_n$ of $A$ such that $A_n \subset D(\delta)$ and
$\bigcup_nA_n$ is dense in $A$.
Furthermore, there exists a sequence $(h_n)$ of elements $h_n$  of the self-adjoint part  $A_{sa}$ of $A$
such that $\delta|_{A_{n}}=\ad(ih_{n})|_{A_{n}}$ (see \cite{BR1}, Example 3.2.25).
If, in addition, $\bigcup_nA_n$ is a core for $\delta_{\alpha}$, i.e.   $\bigcup_nA_n$ is dense
in $D(\delta)$ in the graph norm, then it follows that $\alpha$ is approximately inner.
In particular one has

\begin{equation}
 \lim_{n\rightarrow\infty}\max_{|t|\leq 1}\|\alpha_t(x)-
 \alpha^{(n)}_t(x)\|=0\label{eint1}
\end{equation}
for all $x\in A$ where $\alpha^{(n)}_t=e^{t\delta_n}$ with $\delta_n=\ad(ih_{n})$.
This is a consequence of two general results.
First the strong convergence (\ref{eint1}) of the semigroups is equivalent to strong convergence of
the resolvents $(\iota\pm \delta_n)^{-1}$ to  $(\iota \pm \delta)^{-1}$ by the Kato--Trotter theorem
(see, for example, \cite{Kat1}, Theorem~IX.2.16).
Explicitly (\ref{eint1}) is equivalent to the condition

\begin{equation}
 \lim_{n\to\infty}\|(\iota\pm \delta_n)^{-1}(x)-(\iota \pm \delta)^{-1}(x)\|=0
 \label{eint3}
 \end{equation}
  for all $x$ in a norm-dense subspace of $A$.
 Secondly, since
$\bigcup_{n}A_{n} \subset \cd$, where

\begin{equation}
  \cd=\{x\in D(\delta): \lim_{n\rightarrow\infty}\delta_n(x)=\delta(x)\}\;,\label{eint2}
\end{equation}
it follows that $\cd$ is a core of $\delta$.
In particular the subspaces $(\iota\pm \delta)(\cd)$ are norm-dense  in $A$.
Then if $y\in\cd$ and $x_\pm=(\iota \pm \delta)(y)$  one has
\[
 \lim_{n\to\infty}\|(\iota\pm \delta_n)^{-1}(x_\pm)-(\iota \pm \delta)^{-1}(x_\pm)\|
 \leq  \lim_{n\to\infty}\|(\delta_n-\delta)(y)\|=0
 \;.
 \]
Therefore the strong resolvent convergence  (\ref{eint3}) and the equivalent semigroup convergence (\ref{eint1}) is established.
 Thus the crucial feature in this argument is the core property of $\cd$.
 This follows automatically  if
 $\bigcup_nA_n$  is a core.

   Note that it also follows from general theory that the convergence of the semigroups
   (\ref{eint1}) is equivalent to the graph convergence of the generators, i.e.  for each $x\in D(\delta)$
   there is a sequence $(x_n)$ of $x_n\in A$ such that $\|x_n-x\|\to0$ and $\|\delta_n(x_n)-\delta(x)\|\to0$
   as $n\to\infty$ (see, for example, \cite{BR1}, Theorem~3.1.28).
 Secondly, note  that it does not  follow in general   from Sakai's construction  that $\bigcup_nA_n$ can be taken to be a core for $\delta$.
This was an open problem for many years.
The affirmative answer would imply that $\alpha$ is approximately inner, and this was known as
the Powers-Sakai conjecture.
The problem was finally resolved  in the negative in 2000  (see \cite{K03}, Theorem 1.1).
The counterexample is, however,  an AF-algebra which is not UHF and the problem still seems to be open for UHF-algebras.

The purpose of this note is to analyze the convergence  (\ref{eint1}) for flows
 $\alpha^{(n)}_t=e^{it\delta_n}$ and  $ \alpha_t=e^{it\delta}$ on a general $C^*$-algebra $A$
by examining  the structure of the set $\cd$ defined by (\ref{eint2}).
In particular we consider the relation between $\cd$ and various natural cores of $\delta$.

\section{Cores and spectral subspaces}

Let $\alpha$ be a flow on a $C^*$-algebra $A$ with generator $\delta_\alpha$.
Define the spectral subspace $A^\alpha(K)$, for each  closed subset $K$ of $\Ri$, as the Banach
subspace spanned by the $x\in A$  for which the conditions $f\in L_1(\Ri)$ and
$\supp \hat f\cap K=\emptyset$  imply $\alpha_f(x)=0$.
(Here $\hat f$ denotes the Fourier transform of $f$.)
The $\alpha$-spectrum of $x\in A$, denoted by $\Sp_\alpha(x)$, is defined
 to be the smallest closed subset $K$ of $\Ri$ such that $x\in A^\alpha(K)$.
 Elements with compact $\alpha$-spectra are referred to as geometric elements of $\alpha$
 and the  subspace of geometric elements is denoted by $A^\alpha_G$.
 Explicitly, $A^\alpha_G=\bigcup_{n=1}^\infty A^\alpha([-n,n])$.
Note that each geometric element  is automatically an entire analytic element of $\alpha$.
Moreover, $A^\alpha_G$ is a $^*$-subalgebra of $A$ and an  $\alpha$-invariant core of
$\delta_\alpha$.

 The spectrum $\Sp(\alpha)$ of $\alpha$ is defined to  be the smallest closed subset
 $K$ of $\Ri$ such that  $A=A^\alpha(K)$.
 The Connes spectrum $\Ri(\alpha)$ of $\alpha$ is
 defined to be the intersection of $\Sp(\alpha|_B)$ with $B$ all
 non-zero $\alpha$-invariant hereditary C$^*$-algebras of $A$.
 It is  known that $\Ri(\alpha)$ is a closed subgroup of $\Ri$.

 Our first result establishes under quite general conditions that the geometric elements
 of an approximately inner flow are not contained in the subspace $\cd$ of convergence.

\begin{thm}\label{t1}
Let $A$ be a separable $C^*$-algebra and $\alpha$ an approximately
inner flow with generator $\delta_\alpha$. Suppose there exists a
faithful family of $\alpha$-covariant irreducible representations of
$A$ and that $\Ri(\alpha)\neq \{0\}$. Let $(h_n)$ be a sequence of
self-adjoint elements of $A$ such that
\[
\lim_{n\to\infty}\,\max_{|t|\leq 1}\|\alpha_t(x)-\Ad e^{ith_n}(x)\|=0
\]
for all $x\in A$ and let

\begin{equation}
\cd=\{x\in D(\delta_\alpha):\lim_{n\to\infty}\ad ih_n(x)=\delta_\alpha(x)\}
\;\;\;.
\label{e2.1}
\end{equation}

Then
\[
\cd\not\supset A^\alpha_G\;\;\;.
\]
\end{thm}
\proof\   Under the assumptions of the theorem it follows that there
is a faithful family of irreducible representations $(\pi,\ch_\pi)$
of $A$ such that the representation
\[
\overline \pi=\int^{\oplus}_\Ri dt\,\pi\circ\alpha_t
\]
of $A$ on $L_2(\Ri\,;\ch_\pi)$ {\AK is of type I with} centre
$L_\infty(\Ri)$ if $\Ri(\alpha)=\Ri$ and $L_\infty(\Ri/((2\pi
p)^{-1}\Zi)))$ (as a subalgebra of $L_\infty(\Ri)$) if
$\Ri(\alpha)=p\Zi$. {\AK (This is Theorem 1.2 of \cite{K90} when $A$
is prime. See Remark~\ref{rado1.1} when $A$ is not prime.)} Define a
unitary flow $U$ on $L_2(\Ri\,;\ch_\pi)$ by $U_t\xi(s)=\xi(s+t)$.
Then $U_t\overline{\pi}(x)U_t^*=\overline{\pi}({\alpha}_t(x))$. We
denote by $\overline{\alpha}_t$ the weakly continuous flow $t\mapsto
\Ad\,U_t$ on $\overline{\pi}(A)''$.

Since $A$ is separable, there is a countable faithful family of such
irreducible representations. Let $(\pi_i)$ be such a family.
Let $q\in \Ri(\alpha)$ and define $z_i\in \overline{\pi}_i(A)''\cap
\overline{\pi}_i(A)'\subset L_\infty(\Ri)$ by $z_i(t)=e^{2\pi iqt}$.
Note that $z_i$ satisfies $\overline{\alpha}_t(z_i)=e^{2\pi iqt}z_i$.
Then there is a net $(y_\mu)$ in the unit ball of $A$ such
that $\overline{\pi}_i(y_\mu)$ converges to $z_i$ in the $*$-strong
topology for any $i$ and $\pi(y_\mu)$ converges to $0$ weakly for
any representation $\pi$ disjoint from all $\overline\pi_i$.
We may further suppose that the $\alpha$-spectrum of $y_\mu$ decreases to $\{q\}$.
Since $A$ is separable and the direct sum of
$\overline\pi_i$ is a representation on a separable Hilbert space,
one can choose a sequence $(y_k)$ from the convex combinations of
$(y_\mu)$ such that $(y_k)$ is a central sequence, $\overline
\pi_i(y_k)$ converges to $z_i$ in the $*$-strong topology, and the
$\alpha$-spectrum of $y_k$ shrinks to $\{q\}$.
Since $\overline{\pi}_i(y_ky_k^*)$ converges to the identity for any $i$ and the
direct sum of $\overline{\pi}_i$ is faithful, one may  conclude that
$\|xy_k\|$ converges to $\|x\|$ for any $x\in A$ as
$k\rightarrow\infty$.
We use this fact below.

Let $\varepsilon>0$ and assume $\cd\supset
A^\alpha([-\varepsilon,\varepsilon])$.
 Then it follows from the uniform boundedness theorem that there is constant $c>0$ such that
 \[
 \|\ad ih_m|_{A^\alpha([-\varepsilon,\varepsilon])}\|<c
 \]
 for all $m$. Choose $x\in  A^\alpha([q-\varepsilon/2,q+\varepsilon/2])$ such that $\|x\|=1$
 where $q\in \Ri(\alpha)$ satisfies $q>c+\varepsilon$.
 By the arguments in the previous paragraph one can find a central sequence
 $(y_k)$ in $A$ such that $\|y_k\|\leq1$,
 $\Sp_\alpha(y_k)\subset \langle -q-\varepsilon/2,-q+\varepsilon/2\rangle$ and
 ${\overline \pi}(y_ky_k^*)$ converges weakly to the identity as $k\to\infty$.
 Then $\|xy_k\|\leq 1$, $\|xy_k\|\to1$ and $\Sp_\alpha(xy_k)\subset\langle-\varepsilon,\varepsilon\rangle$.
 Hence $\|\ad ih_m(xy_k)\|<c$ and
 \[
\lim_{k\to\infty} \|\ad ih_m(xy_k)-\ad ih_m(x)y_k\|=0
\;\;\;.
\]
Since $\lim_{m\to\infty}\|\ad ih_m(x)-\delta_\alpha(x)\|=0$ and
$\lim_{k\to\infty}\|\delta_\alpha(x)y_k\|=\|\delta_\alpha(x)\|$ it
follows that $\|\ad ih_m(xy_k)\|>
\|\delta_\alpha(x)\|-\varepsilon/2\geq q-\varepsilon>c$ for all
sufficiently large $m$ and $k$. 
This contradicts the bound $\|\ad
ih_m(xy_k)\|< c$. Therefore
$\cd\not\supset\bigcup^\infty_{n=1}A^\alpha([-n,n])$.
  \hfill$\Box$

\medskip

\begin{remarkn}\label{rado1.7}The foregoing  proof establishes a slightly stronger
 statement: If $\cd$ contains $A^\alpha([-\varepsilon, \varepsilon])$
for some $\varepsilon>0$, then $\cd\cap
A^\alpha([q-\varepsilon/2,q+\varepsilon/2])=\{0\}$ for all large $q\in
\Ri(\alpha)$.
 \end{remarkn}

\medskip

{\AK 
\begin{remarkn}\label{rado1.1} In Theorem~\ref{t1}
we assumed the condition 

\smallskip

 (i) there exists a faithful family of $\alpha$-covariant
irreducible representations of $A$. 

\smallskip

The essential requirement is, however, a consequence of the assumption

\smallskip

 (ii) there exists a faithful family
$\{\pi_i\}$ of irreducible representations of $A$ such that
$\bar{\pi}_i$ is of type I and  the spectrum of $\bar{\alpha}$ on
the center $\bar{\pi}_i(A)''\cap \bar{\pi}_i(A)'$, which we 
denote by $\Delta(\pi)$, is $\Ri(\alpha)$, where $\bar{\pi}_i$ and
$\bar{\alpha}$ are defined in the above proof. 

\smallskip

As we  asserted above (i) implies (ii). 
To confirm this assertion  let us define $\Ri_2(\alpha)$
to be the set of $p\in \Ri$ satisfying: for any non-zero $x\in A$
and any $\varepsilon>0$ there exists an $a\in A^\alpha([p-\varepsilon,
p+\varepsilon])$ such that $\|a\|=1$ and $\|x(a+a^*)x^*\|\geq(2-\varepsilon)\|x\|^2$. 
It is obvious that $\Ri_2(\alpha)$ is a closed subset of $\Ri(\alpha)$. 
Furthermore one can show that the inclusion
$\Ri_2(\alpha)\supset \bigcap_{\pi\in\cf}\Delta(\pi)$ holds for any
faithful family $\cf$ of irreducible representations of $A$ and the
equality holds for some (see   \cite{K87-2}, Proposition 1). 
Consider the conditions

\smallskip

 (i$^\prime$)  $\Ri_2(\hat{\alpha})=\Ri$
 
 \smallskip
 
 and 
 
 \smallskip
 
(ii$^\prime$) $\Ri_2(\alpha)=\Ri$.

\smallskip 

One can show that (i) $\Leftrightarrow$  (i$^\prime$) 
$\Rightarrow$ (ii$^\prime$) $\Leftrightarrow$ (ii).
If $A$ is prime all these conditions are equivalent \cite{K90}. 
The equivalences  of (i) with (i$^\prime$) and (ii) with (ii$^\prime$) are a kind of duality
and these equivalences are straightforward. 
The only implication which is not explicitly
given in the non-prime case seems to be (i$^\prime$) $\Rightarrow$ (ii$^\prime$). 
The arguments we adopt here are given in the proof of Theorem 3.3 of
\cite{K87-1}. 
Let $B=A\times_\alpha\Ri$, let $H$ be a discrete subgroup of
$\Ri(\alpha)$, and let $\beta=\hat{\alpha}|_H$. 
Then it follows from (i') that there is a faithful family of covariant irreducible representations
for $(B\times_\beta H,\hat{H}, \hat{\beta})$ and it follows from $H\subset \Ri(\alpha)$ that 
$H(\hat{\beta})=H$.
Let $x\in B\times_\beta H$. 
For $p\in H^\perp$ and any compact neighborhood $U$ of $p$ in $\Ri$ one
can show from (i$^\prime$) that 
\[
 \sup\{\|x(a+a^*)x^*\|\ :\ a\in B^{\hat{\alpha}}(U),\ \|a\|=1\}=2\,\|x\|^2\;\;\;.
 \]
 Moreover, for $s\in H$ one can show by using Glimm's type of  theorem for the
compact dynamical system $(B\times_\beta H, \hat{H},\hat{\beta})$ \cite{BKR}
that
\[
 \sup\{\|x(a+a^*)x^*\|\ :\ a\in B\lambda(s)\}=2\,\|x\|^2\;\;\;,
\]
where $H\ni s\rightarrow \lambda(s)$ is the canonical unitary group in the multiplier
algebra of $B\times_\beta H$ implementing $\beta$.
Using these two conditions one can  construct a faithful
family $\{\pi_i\}$ of irreducible representations of $B\times_\beta
H$ such that $\pi_i$ restricts to an irreducible representation
$\rho_i$ of $B$ and $\bar{\rho}_i$ which is the direct integral of
$\rho_i\hat{\alpha}_p,\ p\in\Ri$ is of type I with
$\Delta(\rho_i)=H^\perp$. 
Then, by Lemma 5 of \cite{K87-2}, the duality implies that $\Ri_2(\alpha)\supset H$. 
Since $H$ is an arbitrary discrete subgroup of $\Ri(\alpha)$ one can conclude that
$\Ri_2(\alpha)=\Ri(\alpha)$.
 \end{remarkn}
 }

 \medskip

It follows automatically from  the assumptions of Theorem~\ref{t1} that one has
$\cd\not\supset D$ for any core $D$ of $\delta_\alpha$ which contains
the geometric elements $A^\alpha_G$.
In particular $\cd$ cannot contain the analytic elements, or the $C^\infty$ elements, of $\alpha$.
In addition one cannot construct a core in $\cd$ by regularization of the subspace of geometric elements
since the following lemma establishes that the subspace is unchanged
by regularization with respect to the flow.
\begin{lemma}\label{ige}
Let $A^\alpha_G$ denote the geometric elements of the flow $\alpha$ on the $C^*$-algebra $A$.
Then
\begin{eqnarray*}
A^\alpha_G=\{\alpha_f(x): x\in A^\alpha_G\,,\,f\in L_1(\Ri)\,\}
                     =\spann\{\alpha_f(x): x\in A^\alpha_G\,,\,f\in C_c^\infty(\Ri)\,\}
 \;\;\;.
 \end{eqnarray*}
\end{lemma}
\proof\
First $x\in A^\alpha(K)$ if and only if  $g\in L_1(\Ri)$ and
$(\supp \hat g)\cap K=\emptyset$  imply $\alpha_g(x)=0$.
But if $x\in A^\alpha(K)$ and $f\in L_1(\Ri)$ then $\alpha_g(\alpha_f(x))=\alpha_{f*g}(x)$.
Moreover,  $\supp(\widehat{f*g})\subseteq \supp \hat g$.
Therefore  $\alpha_g(\alpha_f(x))=0$ and one deduces that $\alpha_f(x)\in A^\alpha(K)$.
Hence
\[
\{\alpha_f(x): x\in A^\alpha_G\,,\,f\in C_c^\infty(\Ri)\,\}
\subseteq
\{\alpha_f(x): x\in A^\alpha_G\,,\,f\in L_1(\Ri)\,\}
\subseteq A^\alpha_G
\;\;\;.
\]

Next if $x\in A^\alpha(K)$ and $f,g\in L_1(\Ri)$ with $\hat f=\hat g=1$ on an open neighbourhood
$U$ of $K$  then $\alpha_f(x)=\alpha_g(x)$ (see, for example, \cite{BR1}, Lemma~3.2.38).
Replacing $g$ by an approximate identity, with the Fourier transform equal to one on $U$,
and taking the limit one deduces that $\alpha_f(x)=x$.
Therefore
\[
A^\alpha_G=\{\alpha_f(x): x\in A^\alpha_G\,,\,f\in L_1(\Ri)\,\}
\;\;\;.
\]
But the same argument also gives
\[
A^\alpha_G=\{\alpha_f(x): x\in A^\alpha_G\,,\,f\in \cs(\Ri)\,\}
\]
where $\cs(\Ri)$ is the usual Schwartz space.
Finally, if $f\in \cs(\Ri)$ then it follows from the Rubel--Squires--Taylor factorization theorem \cite{RST} (see also \cite{DM})
that there exists a finite set of  $g_i\in C_c^\infty(\Ri)$ and $h_i\in\cs(\Ri)$ such that $f=\sum g_i*h_i$.
In particular $\alpha_f(x)=\sum \alpha_{g_i}(y_i)$ with $y_i=\alpha_{h_i}(x)$.
Therefore
\begin{eqnarray*}
A^\alpha_G&=&\{\alpha_f(x): x\in A^\alpha_G\,,\,f\in \cs(\Ri)\,\}\\[5pt]
&\subseteq &\spann\{\alpha_g(y): y\in A^\alpha_G\,,\,g\in C_c^\infty(\Ri)\,\}
\subseteq \spann\{\alpha_g(y): y\in A^\alpha_G\,,\,g\in \cs(\Ri)\,\}=A^\alpha_G
\end{eqnarray*}
and the proof is complete.
\hfill$\Box$

\bigskip

Despite these observations Example~\ref{ex1} illustrates that many approximately inner flows of interest in  mathematical physics are such that the subspace of convergence $\cd$ contains a dense invariant
set of analytic elements.

\smallskip

It is remarkable that  under the conditions of Theorem~\ref{t1}
convergence of a sequence of bounded derivations  on the subspace
$A^\alpha_G$ automatically implies boundedness of the limit
derivation, at least if $A$ is prime:

\begin{cor} \label{c1}
Let $A$ be a separable prime $C^*$-algebra and $\alpha$ a flow.
Suppose there exists a faithful family of $\alpha$-covariant
irreducible representations of $A$ and that $\Ri(\alpha)\neq \{0\}$.

Suppose that there is a sequence $(b_n)$ of self-adjoint elements of $A$
such that the limits
\[
\delta(x)=\lim_{n\to\infty}\,\ad ib_n(x)
\]
exist for all $x\in A^\alpha_G$.
Then $x\in A^\alpha_G\mapsto\delta(x)\in A$  extends to a  bounded $^*$-derivation on $A$.
\end{cor}
\proof\
It follows from the uniform boundedness theorem that the norm of  $\ad
b_n|_{ A^\alpha([-k,k])}$ is bounded as $n\ra\infty$ for each $k$.
Hence $\delta|_{A^\alpha([-k,k])}$ is bounded.
Then, by \cite{K84}, $\delta$ is closable and its closure $\overline{\delta}$ generates a flow
$\beta$.

By \cite{K02} there is a faithful covariant representation $(\pi,U)$
of $A$ such that the flow $\overline{\alpha}:t\mapsto\Ad U_t$ on the
factor $M=\pi(A)''$ has Connes spectrum $\Ri(\alpha)$, e.g. $\pi$
may be a type II$_\infty$ representation extending the tracial
representation of a UHF algebra (with a non-trivial UHF flow)
``embedded'' in $A${\AK; here we need  $A$ to be prime}. Let
$f\in L_1(\Ri)$ be an integrable real-valued function such that the
Fourier transform $\hat{f}$ has compact support and define
\[
 \delta_f=\int_{\Ri} dt\, f(t)\alpha_t\delta\alpha_{-t}\;.
 \]
The  closure ${\overline\delta}_f$ of $\delta_f$  is also a generator and we denote by $\beta_f$ the flow
generated by ${\overline\delta}_f$.
 By \cite{BK85} the $\alpha$-covariant
representation $\pi$ is also $\beta_f$-covariant. If $f\in L_1(\Ri)$
satisfies $\hat{f}(0)=0$ then $\delta_f$ is bounded.
Moreover,
there is a constant $c>0$ such that $\|\delta_f\|\leq c\|f\|_1$ for
such~$f$.

We fix a positive function $f\in L_1(\Ri)$ such that
$\supp(\hat{f})$ is compact and $\hat{f}(0)=\int dt\, f(t)=1$ and
define $\delta_n$ as the closure of $1/n\int
dt\,f(t/n)\alpha_t\delta\alpha_{-t}$. 
Define $\Delta_n$ to be the
weak extension of $x\mapsto \pi\delta_n(x)$ on $\pi(\cd(\delta_n))$,
i.e.  the generator of the weak extension of the flow on $\pi(A)$
induced by $\delta_n$. 
Since $\|\Delta_n-\Delta_1\|\leq 2c$, we take
a limit point $d$ of the sequence of derivations
$\Delta_n-\Delta_1:\pi(A)''\rightarrow \pi(A)''$ with {\AK pointwise
weak} topology. 
Thus the limit $d$ is a derivation on $M$.
Then by general theory $d$ is an inner derivation on $M$.
Hence $\Delta_1+d$ is a generator. 
Since $\Delta_1+d$ is a limit point of
$\Delta_n:\cd(\Delta_1)\rightarrow M$ and
$\|\overline\alpha_t\Delta_n(x)-\Delta_n\overline\alpha_t(x)\|\ra0,\
x\in M^{\overline\alpha}_G$ as $n\ra\infty$, we conclude that
$\Delta_1+d$ commutes with $\overline{\alpha}$. 
By the arguments originating  in \cite{KR85} we argue that $\Delta_1+d$ generates a
flow which is a bounded perturbation of a scaled $\overline\alpha$
as follows: Let $\gamma$ be the flow on $M$ generated by $\Delta_1+d$. 
We denote by $\overline\alpha\otimes\gamma$ the action
of $\Ri^2$ on $M$ defined by
$(s,t)\mapsto\overline\alpha_s\gamma_t$.
 Since the spectrum of $\overline\alpha\otimes\gamma$ is bounded on each
 $\{p\}\times\Ri$ and the Connes spectrum $\Ri(\overline\alpha\otimes\gamma)$
 is included in $\Ri(\alpha)\times \Ri$, there is a constant $\lambda\in\Ri$ such that
$\Ri(\overline\alpha\otimes\gamma)=\{(p,\lambda p): \, p\in \Ri(\alpha)\}$.
Since $\Sp(\overline\alpha\otimes\gamma)+\Ri(\overline\alpha\otimes\gamma)
=\Sp(\overline\alpha\otimes\gamma)$, we conclude that $\gamma$ is a
bounded perturbation of the flow $t\mapsto \overline\alpha_{\lambda t}$ or the
flow on $A$ generated by $\delta_1$ is a bounded perturbation of the
flow $t\mapsto \alpha_{\lambda t}$.

On the other hand define $\delta_n'$ to be the closure of $n\int
dt\,f(nt)\alpha_t\delta\alpha_{-t}$ and $\Delta'_n$ to be the weak
extension of $x\mapsto \pi(\delta'_n(x))$. 
Then by the same token we
have that $\|\Delta_1-\Delta'_n\|\leq 2c$. 
Note that
$\cd(\Delta_1)=\cd(\Delta'_n)$ and
$\|\pi(\delta_1(x)-\delta'_n(x))\|\leq 2c\|\pi(x)\|$ for $x\in
\cd(\delta_1)$. 
We then conclude that there is a derivation $d'$ of
$M$ such that $\pi(\delta_1(x))-d'(\pi(x))=\pi(\delta(x))$ for $x\in
\pi(A^\alpha_G)$. 
Since $A^\alpha_G$ is dense in $A$, this implies
that $d'$ leaves $\pi(A)$ invariant. 
Hence, since $\pi$ is faithful,
$x\in A^\alpha_G\mapsto \delta_1(x)-\delta(x)$ extends to a
derivation of $A$. 
Combining this with the result in the previous
paragraph, the flow $\beta$ (generated by $\overline\delta$) is a
bounded perturbation of $t\mapsto\alpha_{\lambda t}$. 
If $\lambda\not=0$, this would imply that $\delta_\alpha|_{A^\alpha_G}$
is approximated by inner derivations, which contradicts Theorem~\ref{t1}.
 Hence $\lambda=0$ and $\beta$ is uniformly continuous, or
$\overline\delta$ is bounded. \hfill$\Box$

\bigskip

Theorem~\ref{t1} can be reformulated in various ways. 
The assumption of the existence of a faithful family of $\alpha$-covariant
irreducible representations is {\AK likely to follow from the
approximate innerness of $\alpha$ alone. For example, if each
non-zero ideal of $A$ has a non-zero projection, this follows
because there exist ground states for a perturbed $\alpha$
restricted to an invariant  unital hereditary C$^*$-subalgebra. 
With different arguments we can show this is also} a consequence of a property of
the ideal structure of the $C^*$-algebra.

\begin{prop}\label{p1}
Let $A$ be a separable $C^*$-algebra and $\alpha$ an approximately inner flow.
Suppose that $A$ has at most countably many ideals.
Then there exists a faithful family of $\alpha$-covariant irreducible representations of $A$.
\end{prop}
\proof\ Let $(h_n)$ be a sequence in $A_{sa}$ such that
\[
\alpha_t(x)=\lim_{n\rightarrow\infty}\Ad e^{ith_n}(x)
\]
uniformly in $t$ on every bounded interval of $\Ri$ for all $x\in
A$.

Let $\gamma$ denote the flow on $C_0(\Ri,A)$, the C$^*$-algebra of
continuous functions into $A$ vanishing at infinity,  induced by
translation; $\gamma_t(x)(s)=x(s+t),\ x\in C_0(\Ri,A)$ and let
$\hat{\alpha}$ denote the dual action of $\Ri$ on the crossed
product $A\times_\alpha\Ri$.

 We denote by $\Ni^+$ the one-point compactification of $\Ni$;
$\infty$ is the newly added point.
We assign the dynamical system
$(C_0(\Ri,A),\gamma)$ to each point $n\in\Ni$ and
$(A\times_\alpha\Ri,\hat\alpha)$ to $\infty\in \Ni^+$.
We assert that they form a {\em continuous field of dynamical systems over}
$\Ni^+$.

We define a map $\phi_n$ of $C_c(\Ri,A)$, the space of continuous
functions on $\Ri$ into $A$ with compact support, into $C_0(\Ri,A)$
by
 \[
 \phi_n(f)(p)=\int_{\Ri}dt\, f(t)e^{it(h_n+p)}=\hat{f}(p+h_{n}).
\]
We note that $\phi_n(f)\phi_n(g)=\phi_n(f*_ng)$ and
$\phi_n(f)^*=\phi_n(f^{*_n})$, where
\[
 f*_ng(t)=\int_{\Ri}ds\,  f(s)\Ad e^{ish_n}(g(t-s))
\]
and $f^{*_n}(t)=\Ad e^{ith_n}(f(-t))$.
We denote by $\phi_\infty$ the natural embedding of $C_c(\Ri,A)$ into
$A\times_\alpha\Ri$, which is given by
\[
 \phi_\infty(f)=\int_{\Ri}dt\,  f(t)\lambda_t \;,
\]
where $t\mapsto \lambda_t$ is the canonical unitary flow
implementing $\alpha$ in the multiplier algebra of
$A\times_\alpha\Ri$. We note that
$\phi_\infty(f)\phi_\infty(g)=\phi_\infty(f*g)$ and
$\phi_\infty(f)^*=\phi_\infty(f^*)$, where
\[
 f*g(t)=\int_{\Ri}dt\,  f(s)\alpha_s(g(t-s))
\]
and $f^*(t)=\alpha_t(f(-t))$.
 Note that $f*_ng$ converges to $f*g$
{\AK (resp. $f^{*_n}$ to $f^*$)}  uniformly as continuous functions
of support contained in $\supp(f)+\supp(g)$ {\AK (resp.
$-\supp(f)$). In particular
$\|\phi_n(f*_ng)-\phi_n(f*g)\|\rightarrow0$ and
$\|\phi_n(f^{*_n})-\phi_n(f^*)\|\rightarrow0$.}

For $f\in C_c(\Ri,A)$ and $q\in\Ri$ we define $f_q\in C_c(\Ri,A)$ by
$f_q(t)=f(t)e^{iqt}$. Then $\gamma_q(\phi_n(f))=\phi_n(f_q)$ for
$n\in \Ni$ and $\hat\alpha_q(\phi_\infty(f))=\phi_\infty(f_q)$. The
assertion made above comprises this fact and the following:

\begin{lemma}
If $f\in C_c(\Ri,A)$ then $n\in \Ni^+\mapsto \|\phi_n(f)\|$ is
continuous. The range of $\phi_n$ is dense in $C_0(\Ri,A)$
if $n<\infty$ or in $A\times_\alpha\Ri$ if $n=\infty$.
\end{lemma}

The only non-trivial claim is that
$\lim_{n\rightarrow\infty}\|\phi_n(f)\|=\|\phi_\infty(f)\|$ for
$f\in C_c(\Ri,A)$. Let $\rho(f)=\limsup\|\phi_n(f)\|$. Then it
follows that $\rho$ defines a C$^*$-seminorm on $C_c(\Ri,A)$ as a
$*$-subalgebra of $A\times_\alpha\Ri$. This fact follows from
$\rho(f*g)=\limsup\|\phi_n(f*g)\|=\limsup\|\phi_n(f*_ng)\|\leq
\limsup\|\phi_n(f)\|\|\phi_n(g)\|\leq \rho(f)\rho(g)$,
$\rho(f^**f)=\limsup
\|\phi_n(f^{*_n}*_nf)\|=\limsup\|\phi_n(f)\|^2=\rho(f)^2$ etc.

Since $\rho\hat\alpha_p=\rho$ and $\rho$ is non-zero on a non-zero
element $ag$ with $a\in A$ and $g\in C_c(\Ri)$, one concludes that
$\rho$ is a norm, i.e.  $\rho(f)=\|\phi_\infty(f)\|$. 
Since the same statement holds for any subsequence of $(\phi_n)$, the claim follows.

\medskip

Let $a\in A$ and $g\in C_c(\Ri)$ and define $x\in C_c(\Ri,A)$ by $x(t)=ag(t)$. 
Then $\phi_n(x)(p)=a\int_{\Ri}dt\, g(t)e^{it(h_n+p)}=a\hat g(h_n+p)$. 
Note that $p\mapsto
\|\phi_n(x)(p)\|$ is a continuous function on $\Ri$ {\AK vanishing at $\infty$.}

If $a\not=0$ and $g\not=0$, then $\phi_\infty(x)\not=0$. 
For any $\theta\in \langle 0,\|\phi_\infty(x)\|\rangle$, and for all large
$n$, we find the smallest $p_n\in \Ri$ such that $\|\phi_n(x)(p_n)\|=\theta$. 
Define a seminorm $\rho$ on
$C_c(\Ri,A)$ by
\[
 \rho(f)=\limsup _{n\rightarrow \infty}\max_{p\leq
 p_n}\|\phi_n(f)(p) \|,
\]
which extends to a C$^*$ seminorm on $A\times_\alpha\Ri$. 
Note that $\rho(x)=\theta<\|\phi_\infty(x)\|$ and $q\mapsto
\rho\hat\alpha_q(f)$ is increasing (because
$\gamma_q\phi_n(f)(p)=\phi_n(f_q)=\phi_n(f)(p+q)$).

Let
\[
 \rho_{-\infty}(f)=\lim_{q\rightarrow-\infty}\rho\hat\alpha_q(f),
\]
which defines a C$^*$-seminorm on $A\times_\alpha\Ri$.

Suppose that $\rho_{-\infty}\not=\rho$. Since $\Ker \rho\subsetneqq
\Ker \rho_{-\infty}$, we take an irreducible representation $\pi$ of
the quotient $\Ker\rho_{-\infty}/\Ker\rho$ and regard it as an
irreducible representation of $A\times_\alpha\Ri$.
Then $\Ker\pi\hat\alpha_p\not=\Ker\pi$ for $p\not=0$; if
$\Ker\pi\hat\alpha_p=\Ker\pi$ for some $p\not=0$, then $\Ker\pi$
contains $\Ker\rho\hat\alpha_q$ for any $q$ (as $\Ker\pi\supset
\Ker\rho$), from which follows that $\pi|_{\Ker\rho_{-\infty}}=0$, a
contradiction.
This implies that the center of $\overline\pi$ (as
defined as the direct integral of $\pi\hat\alpha_p$ over $p\in\Ri$
as before) is $L_\infty(\Ri)$, which in turn implies that the
representation $\pi$ of $A\times_\alpha\Ri$ is induced from an
$\alpha$-covariant irreducible representation of $A$.

Suppose that $\rho_{-\infty}=\rho$, which implies that
$\rho\hat\alpha_q=\rho$ for all $q$.
Then there is an
$\alpha$-invariant ideal $I$ of $A$ such that $\Ker\rho$ is
described as $I\times_\alpha\Ri$. For each $\theta\in \langle
0,\|\phi_\infty(x)\|\rangle$ we have defined a seminorm
$\rho=\rho_\theta$ and an ideal $I=I_\theta$ of $A$ if $\rho_\theta$
is $\hat\alpha$-invariant.
Note that if $\rho_\theta$ and
$\rho_{\mu}$ are $\hat\alpha$-invariant for $\theta\not=\mu$, then
$I_\theta\not =I_\mu$. If all $\rho_\theta$ are
$\hat\alpha$-invariant, we thus obtain a continuous family of ideals
of $A$, which contradicts the assumption.
Thus there is a $\theta$
such that $\rho_\theta$ is not $\hat\alpha$-invariant.

Thus we obtain an $\alpha$-covariant irreducible representation $\pi$ of $A$.
If $\Ker\pi$ is non-zero, we apply this argument to
$\Ker\pi$ and $\alpha|_{\Ker\pi}$, which is an approximately inner flow
on a C$^*$-algebra with at most countably many ideals.
By induction one concludes that there is a faithful family of $\alpha$-covariant
irreducible representations. \hfill$\Box$

\section{AF and UHF algebras}

In this section we examine some properties of cores of generators for AF and UHF algebras.
First we note that if $\alpha$ is an approximately inner flow on an AF-algebra then one may choose
the sequence $(h_n)$ which defines the generators of the approximating flows in such a way that the
subspace $\cd$ defined by (\ref{eint2}) is dense in $A$.

\begin{prop}\label{2.3}
Let $A$ be an AF algebra and $\alpha$ an approximately inner flow on
$A$ with generator  $\delta_\alpha$.
Then there exists a sequence $(h_n)$ in $A_{sa}$ such that
\[
 \lim_{n\rightarrow\infty}\max_{|t|\leq 1}\|\alpha_t(x)-\Ad
 e^{ith_n}(x)\|=0
\]
for all $x\in A$ and
\[
  \cd=\{x\in D(\delta_\alpha): \lim_{n\rightarrow\infty}\ad
 ih_n(x)=\delta_\alpha(x)\}
\]
is dense in $A$.
\end{prop}
 \begin{remarkn}\label{rado3.1} Even without assuming the existence of
$(h_n)$ it follows from the result of Sakai mentioned in the
introduction \cite{S76} that there exists an increasing sequence
${A_n}$ of finite dimensional C*-subalgebras of $A$ such that $A_n
\subset D(\delta_\alpha)$ and $\bigcup_nA_n$ is dense in $A$. 
This is true even if $\delta$ is not assumed to be a generator, but only
a closed derivation  (\cite{BR75}, Theorem~11). 
Furthermore, there exists a sequence ${h_n}$ in $A_{sa}$ such that
$\delta|_{A_n}=ad(ih_{n})|_{A_n}$ (see \cite{BR1}, Example 3.2.25).
If $\delta=\delta_{\alpha}$ is a generator and $\bigcup_nA_n$ is a
core for $\delta_{\alpha}$ then the assertion in
Proposition~\ref{2.3} follows with $(h_n)$ equal to this sequence
since in this case $\bigcup_{n}A_{n} \subset \cd$. 
In the absence of the assumption on the existence of $h_{n}$, however, it does not in
general follow from the generator property of $\delta$ that
$\bigcup_nA_n$ can be taken to be a core for $\delta$ or that
${h_{n}}$ exists (see again \cite{K03}, Theorem 1.1).
\end{remarkn}

\noindent {\bf  Proof of Proposition \ref{2.3}.}
 Since $\alpha$ is approximately inner there exists a sequence $d_{n}={d^{*}_{n}}$ in
 $A_{sa}$ such that $\delta_{\alpha}$ is the graph limit of $\ad(id_{n})$.
 This means that for each $x$ in $D(\delta_\alpha)$ there exists a sequence
 $(x_{n})$ of $x_n\in A$ such that

 \begin{equation}
 \lim_{n\to\infty}\|x_{n} - x\| + \lim_{n\to\infty}\|\ad(id_{n})(x_{n})-\delta_{\alpha}(x)\|=0\label{graph}
 \end{equation}
 This condition is actually equivalent to each of the  equivalent conditions
 (\ref{eint1}) and (\ref{eint3}) in the introduction (see \cite{BR1},
 Theorem 3.1.28).

 To prove the proposition we next choose an increasing sequence ${B_n}$ of finite
 dimensional sub-algebras by Sakai's result as in the preceding remark
 such that $\bigcup_{n}B_{n}\subset D(\delta_\alpha)$.
Then we are going to modify $d_{n}$ to $h_{n}$ {\AK by passing to a
subsequence of $d_n$} such that the conclusion of Proposition~\ref{2.3} is valid.

 First, fix an $n$ and let ${\bf X}$ be a set of matrix units for ${B_n}$ (see e.g. \cite{D96},
 Chapter III).
 Now for each matrix unit $x$ in ${B_n}$ there exists, as remarked above,
 a sequence $x_{m}$ in $A$ such that

  \begin{equation}
 \lim_{m\to\infty}\|x_{m} - x\| + \lim_{m\to\infty}\|\ad(id_{m})(x_{m})-\delta_{\alpha}(x)\|=0\label{graph1}
 \end{equation}
  Using Glimm's technique we may furthermore use spectral theory to modify the approximants
 $x_{m}$ for each fixed $m$ such that the set of these approximants form a set of
 matrix units isomorphic to ${\bf X}$ for each $m$. This is possible if we go so far
 out in the sequences indexed by $m$ that the approximation to $x$ is good enough
 (see, for example,  \cite{D96}, Section III.3, or \cite{B72}, Section 2).
 If ${\bf X}_{m}$ denotes the corresponding
 sequences of matrix units, we then have

  \begin{equation}
 \lim_{m\to\infty}||{\bf X}_{m} - {\bf X}||=0
 \end{equation}
 However by careful scrutiny of the proof in \cite{D96} or \cite{B72} one may choose
 the approximants to ${\bf X}$ such that also

  \begin{equation}
 \lim_{m\to\infty}|| \ad (id_{m})({\bf  X}_{m})- \delta_{\alpha}({\bf X})||=0\;\;\;.
 \end{equation}
Hence

 \begin{equation}
 \lim_{m\to\infty}||{\bf X}_{m} - {\bf X}|| + \lim_{m\to\infty}|| \ad (id_{m})({\bf  X}_{m})-
 \delta_{\alpha}({\bf X})||=0
 \label{graph2}
\end{equation}
(The norms in these relations are  defined by taking the maximum of the finite number of norms
 obtained by replacing the sets $\bf X$ by the individual matrix elements $x$ in  $\bf X$.)

 The reason for the convergence of the derivatives is that the matrix elements are
 obtained by applying functional {\AK calculus} by smooth functions in $D(\delta)$. We use
 \cite{BR1}, Theorem~3.2.32,  (or alternatively \cite{P00}, or \cite{ABF90} and \cite{ABF91}) to
 deduce that if $\delta$ is any closed derivation and $x=x^{*}$ is in $D(\delta)$ then
 $f(x)$ is in $D(\delta)$. Since $f(x)$ only depends on the definition of $f$ on ${\rm Spec}(x)$
 we may assume
\[
|||f||| = (2\pi)^{-1/2}\int_{\Ri} dp\, |{\hat f}(p)||p| < \infty
\]
 where $\hat f$ is the Fourier transform of $f$.
 Then

\[
 f(x) = (2{\pi})^{-1/2}\int_{\Ri} dp\,{\hat f}(p)e^{ipx}
 \]
and

\[
 {\delta}(f(x)) = i(2{\pi})^{-1/2}\int_{\Ri} dp {\hat f}(p)p{\int^1_{0}} dt\, e^{itpx}{\delta}(x)e^{i(1-t)px}\;\;\; .
\]
 Hence

\begin{equation}
\|{\delta}(f(x))\| {\le} |||f|||\, \|{\delta}(x)\| .\label{d1}
\end{equation}

We need the following lemma:

 \begin{lemma}\label{conv}
 Let $\delta_{n}$ be a sequence of closed derivations on a C*-algebra $A$
 with graph limit $\delta$ and
 $(x_{n})$ a sequence of elements of  $A$ with  $\lim_{n\to\infty}\|x_{n}-x\|+\|\delta_{n}(x_{n})-\delta(x)\|=0$.
Further let $f_{n}$ be a sequence of functions  converging to $f$ with respect to the semi-norm $|||\cdot|||$.

 It follows that

 \begin{equation}
 \lim_{n\to\infty}\|  \delta_{n}(f_{n}(x_{n})) -  \delta(f(x)) \|=0
 \end{equation}\label{d2}
 \end{lemma}
\proof\
 We have
 \begin{eqnarray*}
 {\delta}_{n}(f_{n}(x_{n})) - {\delta}(f(x))
&=&i(2{\pi})^{-1/2}\int_{\Ri} dp \,{\hat f}_{n}(p)p{\int^1_{0}} dt\, e^{itpx}{\delta}_{n}(x_{n})e^{i(1-t)px}\\
&&\hspace{1cm}-i(2{\pi})^{-1/2}\int_{\Ri} dp \,{\hat f}(p)p{\int^1_{0}} dt\, e^{itpx}{\delta}(x)e^{i(1-t)px}\\
&=& i(2{\pi})^{-1/2}\int_{\Ri} dp\, ({\hat f}_{n}(p) - {\hat f}(p))p{\int^1_{0}} dt\, e^{itpx}{\delta}_{n}(x_{n})e^{i(1-t)px}\\
 &&\hspace{1cm}+i(2{\pi})^{-1/2}\int_{\Ri} dp\, {\hat f}(p)p{\int^1_{0}} dt \,e^{itpx}({\delta_{n}}(x_{n}) - {\delta}(x)) e^{i(1-t)px}
  \end{eqnarray*}
from which one deduces that

\begin{equation}
\|{\delta}_{n}(f_{n}(x_{n})) - {\delta}(f(x))\| \le |||f_{n} - f|||\,\|{\delta}_{n}(x_{n})\|
+|||f|||\, \|\delta_{n}(x_{n})-\delta(x)\|\;\;\;.\label{d3}
\end{equation}
The conclusion of the lemma  follows immediately.\hfill$\Box$

\bigskip

\begin{remarkn}\label{rado3.2}
The semi-norm $|||f|||$ occurring in (\ref{d1}) and (\ref{d3}) can be estimated
by  noting that $|p|(1+p^2)^{1/2}\leq 2\,(1+p^4)^{1/2}$ and using the Cauchy--Schwarz inequality.
Hence

\begin {equation}
|||f|||^2\leq4 \bigg(\int_{\Ri} dp\,(1+p^2)^{-1}\bigg)\,\bigg(\int_{\Ri} dp\,|\hat f(p)|^2(1+p^4)\bigg)
=2\pi\,(\|f''\|_2^2+\|f\|_2^2)
\;\;\;.
\label{d4}
\end{equation}
 It then follows  that the space of functions with $|||f|||<\infty$ contains the Sobolev space
 $W^{2,2}(\Ri)$.
But these estimates are not optimal (see \cite{ABF90},
 \cite{ABF91}, \cite{APS05}, \cite{AZ06}).
 \end{remarkn}

\smallskip

  \noindent {\bf  Proof of Proposition \ref{2.3} continued}$\;$
{\AK We fix a set ${\bf X}^n$ of matrix units for $B_n$ and suppose
that we have obtained a sequence $({\bf X}_{m}^n)$ of matrix units
for each ${\bf X}^n$ by the arguments illustrated above. We then
find a sequence $(u_{n,m})$ of unitaries in $A$ such that

 \begin{equation}
 \lim_{m\to\infty}\|u_{n,m} - \one\|=0\label{graph3}
 \end{equation}
 and

  \begin{equation}
 u_{n,m}{\bf X}^{n}u_{n,m}^{*} = {\bf X}^{n}_{m}   \label{graph4}
 \end{equation}
\noindent where we again interpret (\ref{graph4}) as the set of
relations obtained by replacing ${\bf X}={\bf X}^n_m,{\bf X}^n$ by
each the matrix elements $x$ in $\bf X$. We can choose a subsequence
$(m(n))$ such that $\|u_{n,m(n)}-\one\|<1/n$ and
$\|(\ad(id_{m(n)})\Ad u_{n,m(n)}-\delta)|B_n\|<1/n$. Then it follows
that $\|\hat{\delta_n}(x)-\delta(x)\|\ra0$ for all $x\in \bigcup_k
B_k$ where

 \begin{equation}
\hat{ \delta}_{n}(x) =
u_{n,m(n)}^{*}\ad(id_{m(n)})(u_{n,m(n)}xu_{n,m(n)}^{*})u_{n,m(n)} ,
\label{graph5}
 \end{equation}

\noindent and hence
 \begin{equation}
 \hat{\delta}_{n} = \ad(iu_{n,m(n)}^*d_{m(n)}u_{n,m(n)}) .
 \end{equation}

Since $\|u_{n,m(n)}-1\|\rightarrow1$, one concludes that
$\hat{\delta}_n$ converges to $\delta$ in the graph norm. This
concludes the proof.}




 \smallskip

\begin{remarkn}\label{rado3.3} As we have already said after Proposition \ref{2.3}, if the increasing
 family of finite dimensional sub-algebras constitute a core for $\delta$, there is nothing
 more to prove. So Proposition \ref{2.3}
 only tells something new when the increasing family is not a core, or
 $D(\delta)$ is not an AF Banach algebra in the graph norm. We do not know in general
 whether we can choose approximating inner derivations converging pointwise on a core
 for $\delta$.
 \end{remarkn}

 \smallskip

 The method of constructing the modified sequence again goes back to Glimm, and is
 expanded in Section ~II.3 in \cite{D96} and in Section~2 in \cite{B72}.
   \hfill$\Box$

\bigskip


Although Theorem~\ref{t1} established under quite general conditions that the convergence
subspace $\cd$ cannot contain the analytic elements of the flow $\alpha$ the next example shows
that there are  many examples in which $\cd$ contains an $\alpha$-invariant dense subspace of analytic elements.
The following example is a flow constructed on a one-dimensional lattice.
In mathematical physics terms $A$ is the algebra of observables of a one-dimensional spin-$1/2$ system.
Note that we consider the one-dimensional case for simplicity.
One can  construct similar examples on higher dimensional lattices by analogous arguments.

\begin{exam}\label{ex1}
Let $A$ denote the UHF-algebra given by the $C^*$-closure of the infinite tensor product $\bigotimes_{n\in \Zi}M_2$ of copies of the $2\times 2$-matrices $M_2$.

The algebra $A$ has a natural quasi-local structure.
Let $A_I=\bigotimes_{i\in I}M_2$ denote the family of local matrix algebras indexed by finite subsets $I=\{i_1,\ldots,i_n\}$
with $i_m\in\Zi$.
Further let $A_{\rm loc}=\bigcup_IA_I$.
If $\sigma$ denotes the shift automorphism on $A$
then $\Zi$ acts on $A$ as a group of shifts (space translations)  $n\in\Zi \mapsto \sigma_n=\sigma^n$ which leaves
$A_{\rm loc}$ invariant.
In particular $\sigma_n(A_I)=A_{I+n}$.

Next we construct a flow corresponding to a finite-range interaction between the spins, i.e. 
an interaction which links close by points of $\Zi$.
Fix $\Phi=\Phi^*\in A_J$ for some finite subset $J$.
Define $H_I=H_I^*$ by
\[
H_I=\sum_{m\in I}\sigma_m(\Phi)
\;\;\;.
\]
 Then  introduce the corresponding inner $^*$-derivations $\delta_I$ by
\[
\delta_I(x)= \ad i H_I(x)
\]
for all $x\in A$.
Further define the  $^*$-derivation $\delta$ by $D(\delta)=A_{\rm loc}$ and
\[
\delta(x)=\lim_{I } \ad i H_I(x)
\]
for $x\in A_{\rm loc}$ where the limit is over an increasing family of  $I$  whose union is $\Zi$.
(For each $x\in A_{\rm loc}$ there is an $I\in\Zi$ such that $\delta(x)= \ad i H_I(x)$ by locality.)

There is a flow $\alpha$ on $A$ given by
\[
\alpha_t(x)=\lim_{I}\Ad e^{itH_I }(x)
\;\;\;.
\]
The norm limit exists by the estimates of \cite{BR2}, Theorem 6.2.4.
Moreover the generator $\delta_\alpha$ of the flow is the norm closure $\overline \delta$ of the derivation $\delta$.
Then  $\delta_\alpha$ is the graph limit of the derivations $\delta_I$ (see \cite{BR1}, Theorem~3.1.28).
Therefore if we define
\[
\cd=\{x\in D(\delta_\alpha):\lim_{I} \delta_I(x)=\delta_\alpha(x)\}
\]
one has
\[
A_{\rm loc}=D(\delta)\subseteq \cd\subseteq D(\delta_\alpha)
\;\;\;.
\]
Finally define $\ca$ as the $^*$-algebra generated by $\{\alpha_t(A_{\rm loc}): t\in\Ri\}$.
Then we argue that
$\ca\subseteq \cd$.
To this end it suffices to show that if $x\in A_J$ for some $J$ and $t\in\Ri$ then $\lim_{I} \ad H_I(\alpha_t(x))$ exists.
But if $I_1\subset I_2$ then
\[
\ad H_{I_1}(\alpha_t(x))-\ad H_{I_2}(\alpha_t(x))=\sum_{p\in I_2\backslash I_1}[\sigma_p(\Phi), \alpha_t(x)]
\;\;\;.
\]
Then, as a consequence of \cite{BR2}, Proposition 6.2.9, there are $a, b, c>0$ such that
\[
\|[\sigma_p(\Phi), \alpha_t(x)]\|\leq a\,\|x\|\,e^{-bp+ct}
\]
uniformly for $p\in\Zi$ and $t\in \Ri$.
Therefore
\[
\|\ad H_{I_1}(\alpha_t(x))-\ad H_{I_2}(\alpha_t(x))\|\leq a\,\|x\|\,e^{ct}\sum_{p\in I_1^{\rm c}}e^{-bp}
\]
for all $t\in\Ri$.
It follows immediately that the limit exists and this suffices to establish the inclusion $\ca\subseteq \cd$.

One concludes that $\ca$ is an $\alpha$-invariant subspace of $D(\delta_\alpha)$.
Therefore it is a core of~$\delta_\alpha$.
But it also follows from \cite{BR2}, Theorem 6.2.4, that each $x\in \ca$ is an
analytic element of~$\delta_\alpha$.
Thus $\cd$ contains an $\alpha$-invariant dense subalgebra of analytic elements.

Since the flow $\alpha$ commutes with the group of translations by
$\Zi$  it follows that  the Connes' spectrum $\Ri(\alpha)$, which is
a subgroup of $\Ri$, must be $\Sp(\alpha)$.
Therefore $\Ri(\alpha)\cong\Zi$ or $\Ri$ if $\alpha$ is non-trivial. 
Both cases can occur with a  suitable choice of $\Phi$.

Although the latter conclusions rely on translation invariance one can construct similar examples
on the half line and the same conclusions are valid.
In particular one can add a bounded  $^*$-derivation $\delta_P$ to $\delta_\alpha$ in such a way that $\alpha$
factors into a product of  flows $\alpha^{(\pm)}$  on the left  and right half lattice $\Zi_\pm$, respectively.
Then the algebra $\ca^+$ generated by $\{\alpha^{(+)}_t(A_I) :t\in\Ri, I\subset \Zi_+\}$
is contained in the set $\cd(\alpha^{(+)})$ corresponding to $\alpha^{(+)}$ and consists of analytic elements for the latter flow.
Finally $\Ri(\alpha^{(+)})=\Ri(\alpha)$ because $\alpha^{(-)}\otimes \alpha^{(+)}$ arises by a bounded
perturbation of the generator of $\alpha$.
\hfill$\Box$
\end{exam}

The next example shows that $\cd$ can be much larger  but then the Connes' spectrum is equal to $\{0\}$.
\begin{exam}\label{ex2}
Let $A$ denote the UHF-algebra given by the $C^*$-closure of the infinite tensor product $\bigotimes_{n\geq 1 }M_2$ of copies of the $2\times 2$-matrices $M_2$, $A_I=\bigotimes_{i\in I}M_2$ the  local matrix algebras and $A_{\rm loc}=\bigcup_IA_I$.

Now let $(\lambda_i)_{i\geq 1}$ be a sequence of positive numbers and define $h_i\in A_{\{i\}}$ by
\[
h_i=\left( \begin{array}{cc}
\lambda_i &  0\\
     0             &  0
\end{array}\right)
\;\;\;.
\]
Set $H_n=\sum_{i=1}^nh_i$ and
$\delta_n(x)=\ad i H_n(x)$ for $x\in A$.
Then define  $\alpha$ on $A$  by
\[
\alpha_t(x)=\lim_{I}\Ad e^{itH_I }(x)
\;\;\;.
\]
Set
\[
\cd=\{x\in D(\delta_\alpha):\lim_{I} \delta_I(x)=\delta_\alpha(x)\}
\;\;\;.
\]
We next argue that if the $\lambda_i$ are chosen to increase
sufficiently fast as $i\to\infty$ then
$D(\delta_\alpha^2)\subset\cd$ but in this case $\Ri(\alpha)=\{0\}$.

Let $a_i,a^*_i\in A_{\{i\}}$ be given by
\[
a_i=\left( \begin{array}{cc}
      0  &  1\\
      0  &  0
\end{array}\right)
\;\;\;\;\;\;\;\;{\rm and}\;\;\;\;\;\;\;\;
a_i^*=\left( \begin{array}{cc}
      0  &  0\\
      1  &  0
\end{array}\right)
\]
and note that $H_n=\sum^n_{i=1}\lambda_i\,a^*_ia_i$.
Set $a(I)=\bigotimes_{i\in I}a_i$ and $a^*(J)=\bigotimes_{j\in J}a^*_j$.

Let $C_2$ be the diagonal matrices of $M_2$ and $C$ the
C$^*$-subalgebra generated by $C_2$ at every point of $\Ni$;
$C=\bigotimes_{i\in \Ni}C_2\subset A=\bigotimes_{i\in \Ni}M_2$. For
a subset $K$ of $\Ni$ let $C_K=\bigotimes_{i\in K}C_2$.

We define an action $\gamma$ of $G=\prod_{n=1}^\infty \Ti$ on $A$ by
\[
 \gamma_z=\bigotimes_{n=1}^\infty \Ad \begin{pmatrix}z_n& 0\\ 0&
 1\end{pmatrix}.
\]
Then the fixed point algebra of $\gamma$ is $C$ while $a(I)$ and
$a^*(J)$ are eigen-operators for finite subsets $I,J$:
$\gamma_z(a(I))=\prod_{n\in I}z_n a(I)$. The spectrum of $\gamma$ is
$\coprod\{-1,0,1\}\subset \hat{G}=\coprod_{n\in \Ni} \Zi$, which we
identify with $\cs=\{(I,J)\in P_f(\Ni)\times P_f(\Ni)|\ I\cap
J=\emptyset\}$, where $P_f(\Ni)$ is the set of finite subsets of
$\Ni$ and $p\in \coprod\{-1,0,1\}$ maps to $(I,J)$ with $I=\{n|\
p_n=-1\}$ and $J=\{n|\ p_n=1\}$. Note that for each $(I,J)\in \cs$
the eigen-space is given by $Ca(I)a^*(J)=C_{I^c\cap J^c}a(I)a^*(J)$.

Then each $x\in A_{loc}$ has a unique representation
\[
x=\sum_{(I,J)\in \cs}x(I\,;J)\,a^*(I)a(J)
\]
with $x(I\,;J)\in C_{I^c\cap J^c}$, where the sum is finite. Let
$z^I=\prod_{n\in I}z_n$ and $\bar{z}^J=\prod_{n\in J}\bar{z}_n$.
Since
 \[
 \int_G z^I\bar{z}^J\gamma_z (x)dz= x(I,J)a^*(I)a(J)
\]
with $dz$ is normalized Haar measure on $G$, one deduces that
$\|x(I\,;J)\|\leq \|x\|$.

Now $x\in D(\delta_\alpha^2)$ by locality,
\begin{equation}
\delta_\alpha(x)=i\sum_{(I,J)\in
\cs}(\lambda(I)-\lambda(J))\,x(I\,;J)\,a^*(I)a(J) \label{ex2.1}
\end{equation}
where $\lambda(I)=\sum_{i\in I}\lambda_i$ and
\[
\delta_\alpha^2(x)=-\sum_{(I,J)\in
\cs}(\lambda(I)-\lambda(J))^2\,x(I\,;J)\,a^*(I)a(J) \;\;\;.
\]
In particular $(\lambda(I)-\lambda(J))^2\,\|x(I\,;J)\|\leq \|\delta_\alpha^2(x)\|$.

Next suppose $\lambda_n\geq 2\,(\lambda_1+\ldots+\lambda_{n-1})+6^n$ for all $n$.
If $n=\max (I\cup J)$  then
\[
|\lambda(I)-\lambda(J)|\geq \lambda_n-(\lambda_1+\ldots+\lambda_{n-1})\geq 6^n
\;\;\;.
\]
Therefore $6^n|\lambda(I)-\lambda(J)|
\,\|x(I\,;J)\|\leq \|\delta_\alpha^2(x)\|$ and
\begin{eqnarray*}
\|\delta_\alpha(x)\|&\leq&\sum_{(I,J)\in \cs}|\lambda(I)-\lambda(J)|\,\|x(I\,;J)\|\\[5pt]
&\leq&\sum_{n\geq1}\sum_{\max (I\cup J)=n}6^{-n} \|\delta_\alpha^2(x)\|
\leq \sum_{n\geq1}2^{-n}\|\delta_\alpha^2(x)\|=\|\delta_\alpha^2(x)\|
\end{eqnarray*}
where we have used $\sum_{\max (I\cup J)=n}1<3^n$.
But $A_{\rm loc}$ is a core of $\delta_\alpha^2$.
The foregoing estimates then establish that the representation (\ref{ex2.1})
extends to all $x\in D(\delta_\alpha^2)$;
the infinite sum in (\ref{ex2.1}) is absolutely convergent.

Then if $k\leq l$ one has $H_l-H_k=\sum^l_{i=k+1}\lambda_i\,a^*_ia_i$ and so
\begin{equation}
\delta_l(x)-\delta_k(x)=i\sum_{I,J}(\lambda_{k,l}(I)-\lambda_{k,l}(J))\,x(I\,;J)\,a^*(I)a(J)
\label{ex2.2}
\end{equation}
for all $x\in D(\delta_\alpha^2)$ with $\lambda_{k,l}(I)=\lambda(I\cap\{k+1,\ldots,l\})$.
In particular the summand is only non-zero if $\max (I\cup J)>k$.
But $\lambda_{k,l}(I)=\lambda_{l}(I)-\lambda_{k}(I)$ with $\lambda_{l}(I)=\lambda(I\cap\{1,\ldots,l\})$.
Now if $n+1\in I\cup J$ then
\[
|\lambda_{n+1}(I)-\lambda_{n+1}(J)|\geq \lambda_{n+1}-(\lambda_1+\ldots+\lambda_n)
\]
and
\[
|\lambda_{n}(I)-\lambda_{n}(J)|\leq \lambda_1+\ldots+\lambda_n
\;\;\;.
\]
Hence
\begin{equation}
|\lambda_{n}(I)-\lambda_{n}(J)|\leq |\lambda_{n+1}(I)-\lambda_{n+1}(J)|
\label{ex2.3}
\;\;\;.
\end{equation}
But if $n+1\not\in I\cup J$ then
$\lambda_{n}(I)-\lambda_{n}(J)=\lambda_{n+1}(I)-\lambda_{n+1}(J)$ so
(\ref{ex2.3}) is valid for all $n$. Then by iteration
\[
|\lambda_{k}(I)-\lambda_{k}(J)|\leq |\lambda_{l}(I)-\lambda_{l}(J)|\leq |\lambda(I)-\lambda(J)|
\;\;\;.
\]
Combining these observations one concludes from (\ref{ex2.2}) that
\[
\|\delta_l(x)-\delta_k(x)\|\leq \sum_{\max (I\cup
J)>k}|\lambda(I)-\lambda(J)|\,\|x(I\,;J)\|\leq
2^{-k}\|\delta_\alpha^2(x)\|
\]
for all $x\in D(\delta_\alpha^2)$.
Therefore $\delta_l(x)\to \delta_\alpha(x)$ as $l\to\infty$ and $\cd\supseteq  D(\delta_\alpha^2)$.

Note that in this example $\Ri(\alpha)=\{0\}$ but $\delta_\alpha$ is not bounded.
\hfill$\Box$

\end{exam}

\section*{Acknowledgements}

This work was carried out whilst the first and second authors were visiting
the Centre for Mathematics and its Applications at the Australian National University, Canberra.

\def\CMP{Commun. Math. Phys.\ }
\def\JFA{J. Funct. Anal.}
\def\YMJ{Yokohama Ma\-th. J.\ }
\def\RMP{Rep. Math. Phys.\ }
\def\TMJ{T\^ohoku Ma\-th. J.\ }
\def\RIMS{Publ. RIMS., Kyoto Univ.\ }
\def\JOT{J. Operator Theory\ }
\def\IJM{Internat. J. Math.}

\end{document}